\documentclass[a4, 10pt]{amsart}
\usepackage{amssymb}
\usepackage{amstext}
\usepackage{amsmath}
\usepackage{amscd}
\usepackage{amsthm}
\usepackage{amsfonts}
\usepackage{enumerate}
\usepackage{graphicx}
\usepackage{latexsym}
\usepackage[all]{xy}

\theoremstyle{plain}
\newtheorem{thm}{Theorem}[section]
\newtheorem*{thm*}{Theorem}
\newtheorem*{cor*}{Corollary}
\newtheorem*{prop*}{Proposition}

\newtheorem{lem}[thm]{Lemma}
\newtheorem{cor}[thm]{Corollary}

\newtheorem*{claim*}{Claim}

\theoremstyle{definition}
\newtheorem{dfn}[thm]{Definition}
\newtheorem{rem}[thm]{Remark}
\newtheorem{ex}[thm]{Example}

\newtheorem*{ac}{Acknowledgments}
\newtheorem*{c1}{Case 1}
\newtheorem*{c2}{Case 2}
\newtheorem*{c3}{Case 3}

\theoremstyle{remark}
\newtheorem*{pf}{{\sl Proof}}
\newtheorem*{pff}{Proof of Theorem \ref{main}}

\def\d{\operatorname{d}}
\def\D{{\mathcal D}^b (\mod )}
\def\dual{\operatorname{D}}
\def\e{\epsilon}
\def\exc{\overline{\mathcal E}}
\def\Exc{{\mathcal E}}
\def\Ext{\operatorname{Ext}_{\R}^{1}}
\def\G{\Gamma}
\def\Ga{\Gamma _1}
\def\Gv{\Gamma _0}
\def\gl{\operatorname{gl}\dim}
\def\H{{\mathcal H}}
\def\Hom{\operatorname{Hom}\,}

\def\l{{\mathcal L}}

\def\P{{\Bbb P}}
\def\r{{\mathcal R}}
\def\R{\Lambda}

\def\s{\sigma}
\def\S{\Sigma}
\def\Z{{\Bbb Z}}

\renewcommand{\mod}{{\rm mod\, \R}}

\newcommand{\Gcenter}[2]{
\dimen0=\ht\strutbox
\advance\dimen0\dp\strutbox
\multiply\dimen0 by#1
\divide\dimen0 by2
\advance\dimen0 by-.5\normalbaselineskip
\raisebox{-\dimen0}[0pt][0pt]{#2}}

\begin{document}


\title[]{Exceptional Sequences over path algebras of type $A_n$ and Non-crossing Spanning Trees}
\author[]{Tokuji Araya}
\address{Tokuyama College of Technology, Gakuendai, Shunan, Yamaguchi, 745-8585, Japan}
\email{araya@tokuyama.ac.jp}
\keywords{exceptional sequence, Dynkin quiver, Auslander-Reiten quiver, non-crossing spanning tree}

\subjclass[2000]{16G70, 05C05}
\begin{abstract}
Exceptional sequences are fundamental to investigate the derived categories of finite dimensional algebras.
The aim of this note is to classify all the complete exceptional sequences over the path algebra of a Dynkin quiver of type $A _n$ in terms of non-crossing spanning trees. 
\end{abstract}
\maketitle


\section{Introduction}

The concept of exceptional sequences was introduced by Gorodentsev and Rudakov \cite{GR} to study exceptional vector bundles on $\P ^2$.
Exceptional sequences are certain variation of tilting objects and useful to study the structure of derived categories.
They have effectively been applied to study the derived categories of algebraic varieties by many authors (e.g. \cite{BK}, \cite{R}).
Also exceptional sequences were studied for the path algebras of acyclic quivers by Crawley-Boevey \cite{C} and Ringel \cite{Rin}, for the weighted projective lines by Meltzer \cite{M} and for Cohen-Macaulay modules over the one dimensional graded Gorenstein rings with a simple singularity by the author \cite{A}.
For example, the transitivity of the braid group action was shown for many cases.
But their concrete description of the exceptional sequences does not seem to be studied yet even for the path algebras of Dynkin quivers.

In this paper, we give a combinatorial description of the exceptional sequences in the derived category of the path algebra $\R$ of type $A_n$.
Let us introduce the notation which is necessary to state our results.
Let $\e =(E_1, E_2, \cdots , E_n)$ and $\e '=(E_1', E_2', \cdots , E_n')$ be complete exceptional sequences in $\D$ (see \ref{exdef} for the definition), we write $\e \sim \e '$ if there exists a permutation $\s$ and integers $l_1, l_2,\cdots , l_n$ such that $E_i \cong \S ^{l_i}E_{\s (i)}'$ for every $i$. 
Here, $\S$ is a shift functor in $\D$.
We denote by $\Exc$ the set of complete exceptional sequences in $\D$, and we put $\exc =\Exc /\sim$.

\begin{thm}\label{1.1}
There exists a bijection between $\exc$ and the set of non-crossing spanning trees on the circle with $(n+1)$-points (defined in 2.3).
\end{thm}

We will explicitly give the bijection in terms of the Auslander-Reiten quivers.
The proof of Theorem \ref{1.1} depends on a traditional technique for combinatorics on Auslander-Reiten quivers introduced by Gabriel \cite{G} and extensively studied by Riedtmann \cite{R} and Wiedemann \cite{W} for classification of the representation-finite selfinjective algebras and Gorenstein orders.
Recently, this kind of combinatorics appeared also in cluster tilting theory by Caldero-Chapoton-Schiffler \cite{CCS} and Iyama \cite{I}.

It is known the number of non-crossing spanning trees with $(n+1)$-points by Dulucq and Penaud \cite{DP} (see also \cite{FN}).
Therefore we get following.

\begin{cor}
The cardinality of $\exc$ is equal to ${\displaystyle \frac{1}{2n+1}{3n \choose n}}$.\end{cor}

On the other hand, Seidel \cite{S} determines that the number of exceptional sequences (in $\mod$) is equal to $(n+1)^{n-1}$ by induction on $n$,
but a concrete description of the exceptional sequences is not given there.

Also we will show in Section 5 that there is a nice interpretation of mutation of exceptional sequences in terms of non-crossing spanning trees.


\section{Preliminaries}

Let $\R$ be the path algebra of a Dynkin quiver of type $A _n$ over a field $k$.
We denote by $\mod$ the category of finitely generated right $\R$-modules and by $\D$ the bounded derived category.
We write $\S$ the shift functor of $\D$.

\begin{dfn}\label{exdef}
An object $E \in \D$ is called {\it exceptional} if $\Hom (E,E) \cong k$ and $\Hom (\S ^lE,E)=0$ for all integers $l \not= 0$.
A pair $(E,F)$ of exceptional objects is called an {\it exceptional pair} if $\Hom (\S ^l F,E)=0$ for all integers $l$.
A sequence $\e =(E_1,E_2, \cdots ,E_r)$ of exceptional objects is called an {\it exceptional sequence} of length $r$ 
if $(E_i, E_j)$ is an exceptional pair for each $i<j$.
An exceptional sequence $\e$ is called {\it complete} if the length of $\e$ is equal to the number $n$ of simple $\R$-modules.
\end{dfn}

The derived category has the following properties.

\begin{rem}\label{2.2}
\begin{enumerate}[1.]
\item The derived category $\D$ is independent of the choice of orientation of $A_n$ \cite{H}.
\item In general, exceptional objects are indecomposable.
If $\R$ is the path algebra of a Dynkin quiver, the converse also holds.
Thus, in our situation that $\R$ is the path algebra of $A_n$, 
an object is exceptional if and only if it is indecomposable.
\item For any indecomposable object $X\in \D$, there exists a unique integer $l$ and a unique indecomposable $\R$-module $Y$ such that $X \cong \S ^lY$ \cite{H}.
\item Let $(E_1,E_2,\cdots,E_n)$ be an exceptional sequence, then so is $(\S^{l_1}E_1,\S^{l_2}E_2,\cdots,$ $\S^{l_n}E_n)$ for any integers $l_1,l_2,\cdots,l_n$.
\item Since $\gl \R \leq 1$, we see that $E \in \mod $ is exceptional if $\Hom (E,E) \cong k$ and $\Ext (E,E)=0$.
Furthermore, for exceptional modules $E$ and $F$, we can see that $(E,F)$ is an exceptional pair if $\Hom (F,E)=\Ext (F,E)=0$.
\end{enumerate}
\end{rem}

By Remark \ref{2.2}.1, we only have to consider the case when the orientation is $$\overset{1}{\bullet} \leftarrow \overset{2}{\bullet} \leftarrow \cdots \leftarrow \overset{n}{\bullet}.$$
For any $(E_1, E_2, \cdots , E_n) \in \Exc$, by Remark \ref{2.2}.3 and 4, we may assume that all $E_i$ are finitely generated $\R$-modules.

Let $\G =(\Gv ,\Ga)$ be the Auslander-Reiten quiver (c.f. \cite{ARS}, \cite{ASS}) of $\mod$.
We identify the set $\Gv$ of vertices in $\G$ with the isomorphism classes of indecomposable $\R$-modules. 
For each integers $0 \leq i < j \leq n$, we write $X_{i,j}$ the indecomposable $\R$-module whose representation is given by $( \overset{1}{0} \leftarrow \cdots \leftarrow \overset{i}{0} \leftarrow \overset{i+1}{k} \leftarrow \cdots \leftarrow \overset{j}{k} \leftarrow \overset{j+1}{0} \leftarrow \cdots \leftarrow \overset{n}{0})$.
Then it is known that $\Gv = \{ X_{i,j} \ | \ 0 \leq i < j \leq n \}$ and $\Gamma$ is the following form:

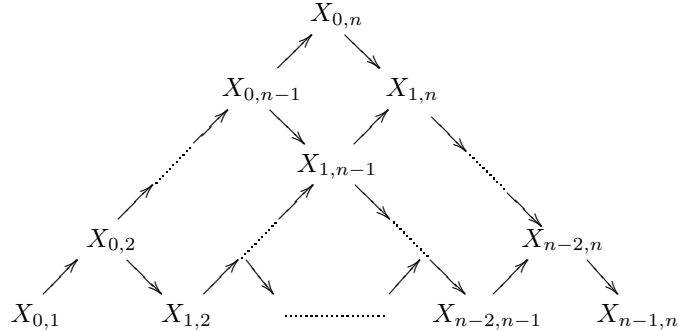
\[
    \begin{xy}
    (0,0)*{X_{0,1}},<1cm,1cm>*{X_{0,2}},<3cm,3cm>*{X_{0,n-1}},<4cm,4cm>*{X_{0,n}},
                   <2cm,0cm>*{X_{1,2}},<4cm,2cm>*{X_{1,n-1}},<5cm,3cm>*{X_{1,n}},
                   <6cm,0cm>*{X_{n-2,n-1}},<7cm,1cm>*{X_{n-2,n}},<8cm,0cm>*{X_{n-1,n}},
{\ar (1,3);(5.5,7.5)}{\ar (11,13);(15.5,17.5)}{\ar @{.} (16,18);(20.5,22.5)}{\ar (21,23);(25.5,27.5)}
{\ar (31,33);(35.5,37.5)}{\ar (11,7.5);(15.5,3)}{\ar (30,27.5);(34.5,23)}{\ar (40,37.5);(44.5,33)}
{\ar (20,3);(24.5,7.5)}{\ar @{.} (25,8);(29.5,12.5)}{\ar (30,13);(34.5,17.5)}{\ar (40,23);(44.5,27.5)}
{\ar (39,17.5);(43.5,13)}{\ar @{.} (44,12.5);(48.5,8)}{\ar (49,7.5);(53.5,3)}
{\ar (48,26.5);(52.5,22)}{\ar @{.} (53,21.5);(57.5,17)}{\ar (58,16.5);(62.5,12)}
{\ar (55,3);(59.5,7.5)}{\ar (67.5,7.5);(71.5,3)}
{\ar (21,7.5);(24.5,3)}{\ar (40,3);(44,7.5)}
{\ar @{.} (25,0);(38,0)}
    \end{xy}
\]

Now we consider the circle with $n+1$ points labeled $0,1, 2,\cdots , n$ counter clockwise on it.
We write $c(i,j) (=c(j,i))$ the chord between the points $i$ and $j$.
We denote by $C_{n+1}$ the set of chords on the circle.

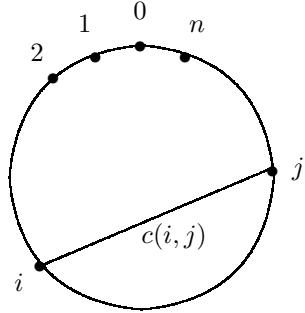
\begin{center}
\begin{picture}(140,140)

\qbezier (10,65) (14,111) (60,115)
\qbezier (10,65) (14,19) (60,15)
\qbezier (60,115) (106,111) (110,65)
\qbezier (60,15) (106,19) (110,65)

\put (57,112){$\bullet$}
\put (40,108){$\bullet$}
\put (74,108){$\bullet$}
\put (24,100){$\bullet$}
\put (19,29){$\bullet$}
\put (107,65){$\bullet$}

\put (57,125){$0$}
\put (36,120){$1$}
\put (18,109){$2$}
\put (12,22){$i$}
\put (117,66){$j$}
\put (78,120){$n$}

\qbezier (23,32) (23,32) (108,68)

\put (60,40){$c(i,j)$}

\end{picture}
\end{center}

Since $\Gv =\{ X_{i,j} | 0 \leq i < j \leq n\}$ and $C_{n+1}=\{ c(i,j) | 0 \leq i < j \leq n\}$, there is a bijection 

$$
\Phi : \Gv \to C_{n+1}\ \ (X_{i,j} \mapsto c(i,j)).
$$

\begin{dfn}
A set $T$ consisting of $n$ chords in $C_{n+1}$ is called a {\it non-crossing spanning tree} if the chords in $T$ form a tree.
\end{dfn}

We remark that each vertex is then contained in at least one chord.

Here we give some examples of non-crossing spanning tree and not non-crossing spanning trees on the circle with $6$ points.

\begin{ex}

The first figure is a non-crossing spanning tree, but last two figures are not non-crossing spanning trees because they have a cycle.

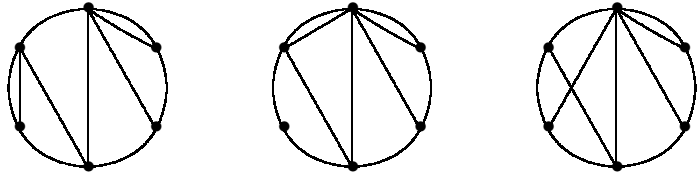
\begin{picture}(270,90)(-20,0)

\qbezier (7,40) (9,68) (37,70)
\qbezier (7,40) (9,12) (37,10)
\qbezier (37,10) (65,12) (67,40)
\qbezier (37,70) (65,68) (67,40)

\put (35,68){$\bullet$}
\put (9,53){$\bullet$}
\put (9,23){$\bullet$}
\put (35,8){$\bullet$}
\put (60.5,23){$\bullet$}
\put (60.5,53){$\bullet$}

\qbezier (11.5,55) (11.5,40) (11.5,25)
\qbezier (11.5,55) (24.25,32.5) (37,10)
\qbezier (37,70) (37,40) (37,10)
\qbezier (37,70) (49.75,47.5) (62.5,25)
\qbezier (37,70) (47.75,62.5) (62.5,55)


\qbezier (107,40) (109,68) (137,70)
\qbezier (107,40) (109,12) (137,10)
\qbezier (137,10) (165,12) (167,40)
\qbezier (137,70) (165,68) (167,40)

\put (135,68){$\bullet$}
\put (109,53){$\bullet$}
\put (109,23){$\bullet$}
\put (135,08){$\bullet$}
\put (160.5,23){$\bullet$}
\put (160.5,53){$\bullet$}

\qbezier (111.5,55) (124.25,62.5) (137,70)
\qbezier (111.5,55) (124.25,32.5) (137,10)
\qbezier (137,70) (137,40) (137,10)
\qbezier (137,70) (149.75,47.5) (162.5,25)
\qbezier (137,70) (147.75,62.5) (162.5,55)


\qbezier (207,40) (209,68) (237,70)
\qbezier (207,40) (209,12) (237,10)
\qbezier (237,10) (265,12) (267,40)
\qbezier (237,70) (265,68) (267,40)

\put (235,68){$\bullet$}
\put (209,53){$\bullet$}
\put (209,23){$\bullet$}
\put (235,8){$\bullet$}
\put (260.5,23){$\bullet$}
\put (260.5,53){$\bullet$}

\qbezier (237,70) (224.25,47.5) (211.5,25)
\qbezier (211.5,55) (224.25,32.5) (237,10)
\qbezier (237,70) (237,40) (237,10)
\qbezier (237,70) (249.75,47.5) (262.5,25)
\qbezier (237,70) (247.75,62.5) (262.5,55)

\end{picture}
\end{ex}

Now we can state Theorem \ref{1.1} in the following more explicit form.

\begin{thm}\label{main}
A bijection from $\exc$ to the set of non-crossing
spanning trees on the circle with $(n+1)$-points is given by
\[\epsilon=(E_1,E_2,\cdots,E_n)\mapsto
\Phi(\epsilon):=\{\Phi(E_i)\ |\ 1\le i\le n\}.\]
\end{thm}

We shall give a proof of Theorem \ref{main} in the next section.


\section{Proof of Theorem \ref{main}}

For $X \in \Gv$, we consider the following four subsets of $\Gv$.
\begin{eqnarray*}
\H^0_{+}(X) & = & \{ Y \in \Gv|\ \Hom (X,Y)\not= 0\}, \\
\H^0_{-}(X) & = & \{ Y \in \Gv|\ \Hom (Y,X)\not= 0\}, \\
\H^1_{+}(X) & = & \{ Y \in \Gv|\ \Ext (X,Y)\not= 0\}, \\
\H^1_{-}(X) & = & \{ Y \in \Gv|\ \Ext (Y,X)\not= 0\},
\end{eqnarray*}
and we put $\H_{\pm}(X) = \H^0_{\pm}(X) \cup \H^1_{\pm}(X)$.
The following lemma is obtained by \cite[Lemma 8.1]{Hi} and by the Auslander-Reiten duality (\cite{ARS}, \cite{ASS}).

\begin{lem}\label{3.1}
For $X=X_{i,j}$, the above four sets are given by the following forms in $\G$.

\begin{eqnarray*}
\H^0 _{+}(X_{i,j}) & = & \{ X_{s,t} |\ i \leq s \leq j-1,\ j \leq t \leq n \}, \\
\H^0 _{-}(X_{i,j}) & = & \{ X_{s,t} |\ 0 \leq s \leq i,\ i+1 \leq t \leq j \}, \\
\H^1 _{+}(X_{i,j}) & = & \{ X_{s,t} |\ 0 \leq s \leq i-1,\ i \leq t \leq j-1 \}, \\
\H^1 _{-}(X_{i,j}) & = & \{ X_{s,t} |\ i+1 \leq s \leq j,\ j+1 \leq t \leq n \}.
\end{eqnarray*}

In particular, we can draw the areas $\H^0_{\pm}(X_{i,j})$ and $\H^1_{\pm}(X_{i,j})$ in $\G$.

{\tiny
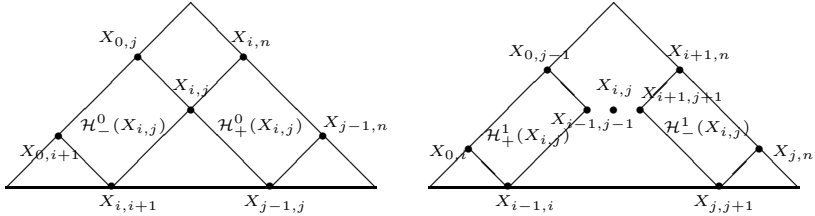
\begin{picture}(300,100)

\put (0,10){\line(1,1){70}}
\put (70,80){\line(1,-1){70}}
\put (0,10){\line(1,0){140}}

\put (50,60){\line(1,-1){50}}
\put (40,10){\line(1,1){50}}
\put (20,30){\line(1,-1){20}}
\put (100,10){\line(1,1){20}}

\put (63,47){$X_{i,j}$}

\put (5,22){$X_{0,i+1}$}
\put (35,3){$X_{i,i+1}$}
\put (35,65){$X_{0,j}$}

\put (85,65){$X_{i,n}$}
\put (90,3){$X_{j-1,j}$}
\put (120,33){$X_{j-1,n}$}

\put (80,32){$\H^0_{+}(X_{i,j})$}
\put (28,32){$\H^0_{-}(X_{i,j})$}

\put (68,38){$\bullet$}

\put (48,58){$\bullet$}
\put (38,9){$\bullet$}
\put (18,28){$\bullet$}

\put (88,58){$\bullet$}
\put (98,9){$\bullet$}
\put (118,28){$\bullet$}

\put (160,10){\line(1,1){70}}
\put (230,80){\line(1,-1){70}}
\put (160,10){\line(1,0){140}}

\put (175,25){\line(1,-1){15}}
\put (190,10){\line(1,1){30}}
\put (205,55){\line(1,-1){15}}
\put (240,40){\line(1,1){15}}
\put (240,40){\line(1,-1){30}}
\put (270,10){\line(1,1){15}}

\put (223,47){$X_{i,j}$}

\put (206,34){$X_{i-1,j-1}$}
\put (160,22){$X_{0,i}$}
\put (185,3){$X_{i-1,i}$}
\put (190,60){$X_{0,j-1}$}
\put (240,44){$X_{i+1,j+1}$}
\put (250,60){$X_{i+1,n}$}
\put (260,3){$X_{j,j+1}$}
\put (290,22){$X_{j,n}$}

\put (249,31){$\H^1_{-}(X_{i,j})$}
\put (181,28){$\H^1_{+}(X_{i,j})$}

\put (228,38){$\bullet$}

\put (218,38){$\bullet$}
\put (203,53){$\bullet$}
\put (188,9){$\bullet$}
\put (173,23){$\bullet$}
\put (238,38){$\bullet$}
\put (253,53){$\bullet$}
\put (268,9){$\bullet$}
\put (283,23){$\bullet$}

\end{picture}
}

\end{lem}

For the points $i$ and $j$ on the circle, we can consider the distance $\d(i,j)$ from $i$ to $j$ obtained by a map $\d :\{ 0,1,\cdots ,n\} \times \{ 0,1,\cdots ,n\} \to {\Bbb N}$ which is given by $$\d (i,j)= \left\{ \begin{array}{ll}
j-i & (j \geq i) \\
j+n+1-i & (j < i).
\end{array}
\right.$$

We can easily check the following lemma by the definition of exceptional pair and by Lemma \ref{3.1}.
It is a key to prove Theorem \ref{main}.

\begin{lem}\label{3.2}
Let $X=X_{i,j}, X'=X_{i',j'} \in \Gv$.
\begin{enumerate}[1.]
\item Both $(X,X')$ and $(X',X)$ are exceptional pairs if and only if $\Phi (X)$ does not meet $\Phi (X')$.
\item Neither $(X,X')$ nor $(X',X)$ is an exceptional pair if and only if $\Phi (X)$ meets $\Phi (X')$ in the interior of the circle.
\item $(X,X')$ is an exceptional pair but $(X',X)$ is not an exceptional pair if and only if one of the following conditions holds:
\begin{enumerate}[(i)]
\item $i=i'$ and $\d (i,j) < \d (i',j')$,
\item $i=j'$ and $\d (i,j) < \d (j',i')$,
\item $j=i'$ and $\d (j,i) < \d (i',j')$,
\item $j=j'$ and $\d (j,i) < \d (j',i')$.
\end{enumerate}
\end{enumerate}
\end{lem}

Now we can prove Theorem \ref{main}.

\begin{pff}
Let $\e=(E_1,E_2,\cdots,E_n)$ be a complete exceptional sequence.
It comes from Definition \ref{exdef} and Lemma \ref{3.2}, one can easily check that  $\Phi (\e )=\{ \Phi (E_1),\Phi (E_2),\cdots \Phi (E_n) \}$ does not have cycles.
Therefore $\Phi (\e )$ is a non-crossing spanning tree.

Let $T=\{ c_1,c_2, \cdots ,c_n \}$ be a non-crossing spanning tree.
Since $c_i$ does not meet $c_j$ in the interior of the circle, either $(\Phi ^{-1}(c_i),\Phi ^{-1}(c_j))$ or $(\Phi ^{-1}(c_j),\Phi ^{-1}(c_i))$ is an exceptional pair for all $i \not= j$.
We remark that there exists $c_{i_1} \in T$ such that $(\Phi ^{-1}(c_{i_1}),\Phi ^{-1}(c))$ is an exceptional pair for all $c \in T \setminus \{ c_{i_1} \}$.
Indeed, if such $c_{i_1}$ does not exist, then $T$ must have a cycle.
We put $E_1=\Phi ^{-1}(c_{i_1})$.
We choose $c_{i_j}$ such that $(\Phi ^{-1}(c_{i_j}),\Phi ^{-1}(c))$ is an exceptional pair for all $c \in T \setminus \{ c_{i_1},c_{i_2},\cdots ,c_{i_j} \}$ inductively and put $E_j=\Phi ^{-1}(c_{i_j})$.
Then $(E_1, E_2, \cdots ,E_n)$ is a complete exceptional sequence.

\qed
\end{pff}

\begin{ex}
If $n=3$, the following quiver is the Auslander-Reiten quiver of $\mod$.

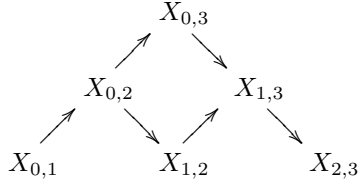
\[
    \begin{xy}
    (0,0)*{X_{0,1}},<1cm,1cm>*{X_{0,2}},<2cm,2cm>*{X_{0,3}},
<2cm,0cm>*{X_{1,2}},<3cm,1cm>*{X_{1,3}},<4cm,0cm>*{X_{2,3}},
{\ar (1,3);(5.5,7.5)}{\ar (11,13);(15.5,17.5)}{\ar (20,3);(24.5,7.5)}
{\ar (11,7.5);(15.5,3)}{\ar (20,17.5);(24.5,13)}{\ar (30,7.5);(34.5,3)}
    \end{xy}
\]

In this case, there are $16$ complete exceptional sequences and $12$ non-crossing spanning trees.
The following are the complete exceptional sequences and the corresponding non-crossing spanning trees.

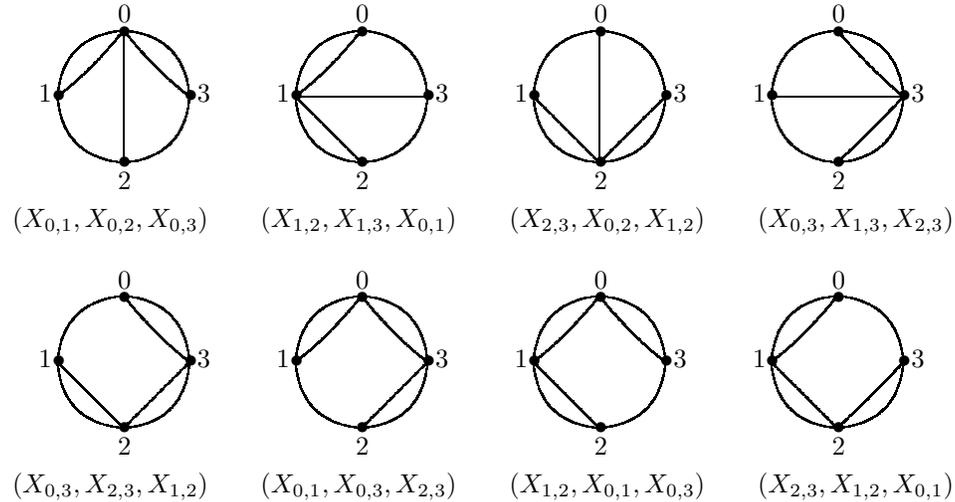
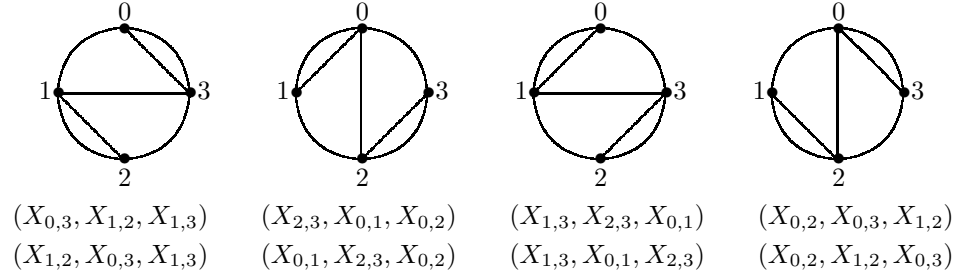
\begin{picture}(340,80)

\qbezier (7,35) (9,58) (32,60)
\qbezier (7,35) (9,12) (32,10)
\qbezier (32,60) (55,58) (57,35)
\qbezier (32,10) (55,12) (57,35)

\put (30,57){$\bullet$}
\put (5,33){$\bullet$}
\put (30,8){$\bullet$}
\put (55,33){$\bullet$}

\put (30,63){$0$}
\put (0,33){$1$}
\put (30,0){$2$}
\put (60,33){$3$}

\qbezier (7,35) (19.5,44.5) (32,60)
\qbezier (32,10) (32,35) (32,60)
\qbezier (32,60) (44.5,44.5) (57,35)

\qbezier (97,35) (99,58) (122,60)
\qbezier (97,35) (99,12) (122,10)
\qbezier (122,60) (145,58) (147,35)
\qbezier (122,10) (145,12) (147,35)

\put (120,57){$\bullet$}
\put (95,33){$\bullet$}
\put (120,8){$\bullet$}
\put (145,33){$\bullet$}

\put (120,63){$0$}
\put (90,33){$1$}
\put (120,0){$2$}
\put (150,33){$3$}

\qbezier (97,35) (109.5,44.5) (122,60)
\qbezier (97,35) (122,35) (147,35)
\qbezier (97,35) (109.5,22.5) (122,10)

\qbezier (187,35) (189,58) (212,60)
\qbezier (187,35) (189,12) (212,10)
\qbezier (212,60) (235,58) (237,35)
\qbezier (212,10) (235,12) (237,35)

\put (210,57){$\bullet$}
\put (185,33){$\bullet$}
\put (210,8){$\bullet$}
\put (235,33){$\bullet$}

\put (210,63){$0$}
\put (180,33){$1$}
\put (210,0){$2$}
\put (240,33){$3$}

\qbezier (187,35) (199.5,22.5) (212,10)
\qbezier (212,10) (212,35) (212,60)
\qbezier (212,10) (224.5,22.5) (237,35)

\qbezier (277,35) (279,58) (302,60)
\qbezier (277,35) (279,12) (302,10)
\qbezier (302,60) (325,58) (327,35)
\qbezier (302,10) (325,12) (327,35)

\put (300,57){$\bullet$}
\put (275,33){$\bullet$}
\put (300,8){$\bullet$}
\put (325,33){$\bullet$}

\put (300,63){$0$}
\put (270,33){$1$}
\put (300,0){$2$}
\put (330,33){$3$}

\qbezier (302,60) (315.5,44.5) (327,35)
\qbezier (277,35) (302,35) (327,35)
\qbezier (302,10) (315.5,22.5) (327,35)

\end{picture}

\begin{center}
$
(X_{0,1},X_{0,2},X_{0,3})\hspace{7mm}(X_{1,2},X_{1,3},X_{0,1})\hspace{7mm}(X_{2,3},X_{0,2},X_{1,2})\hspace{7mm}(X_{0,3},X_{1,3},X_{2,3})$
\end{center}

\begin{picture}(340,80)

\qbezier (7,35) (9,58) (32,60)
\qbezier (7,35) (9,12) (32,10)
\qbezier (32,60) (55,58) (57,35)
\qbezier (32,10) (55,12) (57,35)

\put (30,57){$\bullet$}
\put (5,33){$\bullet$}
\put (30,8){$\bullet$}
\put (55,33){$\bullet$}

\put (30,63){$0$}
\put (0,33){$1$}
\put (30,0){$2$}
\put (60,33){$3$}

\qbezier (7,35) (19.5,22.5) (32,10)
\qbezier (32,10) (44.5,22.5) (57,35)
\qbezier (32,60) (44.5,44.5) (57,35)

\qbezier (97,35) (99,58) (122,60)
\qbezier (97,35) (99,12) (122,10)
\qbezier (122,60) (145,58) (147,35)
\qbezier (122,10) (145,12) (147,35)

\put (120,57){$\bullet$}
\put (95,33){$\bullet$}
\put (120,8){$\bullet$}
\put (145,33){$\bullet$}

\put (120,63){$0$}
\put (90,33){$1$}
\put (120,0){$2$}
\put (150,33){$3$}

\qbezier (97,35) (109.5,44.5) (122,60)
\qbezier (122,60) (134.5,44.5) (147,35)
\qbezier (122,10) (134.5,22.5) (147,35)

\qbezier (187,35) (189,58) (212,60)
\qbezier (187,35) (189,12) (212,10)
\qbezier (212,60) (235,58) (237,35)
\qbezier (212,10) (235,12) (237,35)

\put (210,57){$\bullet$}
\put (185,33){$\bullet$}
\put (210,8){$\bullet$}
\put (235,33){$\bullet$}

\put (210,63){$0$}
\put (180,33){$1$}
\put (210,0){$2$}
\put (240,33){$3$}

\qbezier (187,35) (199.5,22.5) (212,10)
\qbezier (187,35) (199.5,44.5) (212,60)
\qbezier (212,60) (224.5,44.5) (237,35)

\qbezier (277,35) (279,58) (302,60)
\qbezier (277,35) (279,12) (302,10)
\qbezier (302,60) (325,58) (327,35)
\qbezier (302,10) (325,12) (327,35)

\put (300,57){$\bullet$}
\put (275,33){$\bullet$}
\put (300,8){$\bullet$}
\put (325,33){$\bullet$}

\put (300,63){$0$}
\put (270,33){$1$}
\put (300,0){$2$}
\put (330,33){$3$}

\qbezier (277,35) (289.5,44.5) (302,60)
\qbezier (277,35) (289.5,22.5) (302,10)
\qbezier (302,10) (315.5,22.5) (327,35)

\end{picture}

\begin{center}
$
(X_{0,3},X_{2,3},X_{1,2})\hspace{7mm}(X_{0,1},X_{0,3},X_{2,3})\hspace{7mm}(X_{1,2},X_{0,1},X_{0,3})\hspace{7mm}(X_{2,3},X_{1,2},X_{0,1})$
\end{center}

\begin{picture}(340,80)

\qbezier (7,35) (9,58) (32,60)
\qbezier (7,35) (9,12) (32,10)
\qbezier (32,60) (55,58) (57,35)
\qbezier (32,10) (55,12) (57,35)

\put (30,57){$\bullet$}
\put (5,33){$\bullet$}
\put (30,8){$\bullet$}
\put (55,33){$\bullet$}

\put (30,63){$0$}
\put (0,33){$1$}
\put (30,0){$2$}
\put (60,33){$3$}

\qbezier (7,35) (19.5,22.5) (32,10)
\qbezier (7,35) (32,35) (57,35)
\qbezier (32,60) (44.5,47.5) (57,35)

\qbezier (97,35) (99,58) (122,60)
\qbezier (97,35) (99,12) (122,10)
\qbezier (122,60) (145,58) (147,35)
\qbezier (122,10) (145,12) (147,35)

\put (120,57){$\bullet$}
\put (95,33){$\bullet$}
\put (120,8){$\bullet$}
\put (145,33){$\bullet$}

\put (120,63){$0$}
\put (90,33){$1$}
\put (120,0){$2$}
\put (150,33){$3$}

\qbezier (97,35) (109.5,47.5) (122,60)
\qbezier (122,10) (122,35) (122,60)
\qbezier (122,10) (134.5,22.5) (147,35)

\qbezier (187,35) (189,58) (212,60)
\qbezier (187,35) (189,12) (212,10)
\qbezier (212,60) (235,58) (237,35)
\qbezier (212,10) (235,12) (237,35)

\put (210,57){$\bullet$}
\put (185,33){$\bullet$}
\put (210,8){$\bullet$}
\put (235,33){$\bullet$}

\put (210,63){$0$}
\put (180,33){$1$}
\put (210,0){$2$}
\put (240,33){$3$}

\qbezier (187,35) (199.5,47.5) (212,60)
\qbezier (187,35) (212,35) (237,35)
\qbezier (212,10) (224.5,22.5) (237,35)

\qbezier (277,35) (279,58) (302,60)
\qbezier (277,35) (279,12) (302,10)
\qbezier (302,60) (325,58) (327,35)
\qbezier (302,10) (325,12) (327,35)

\put (300,57){$\bullet$}
\put (275,33){$\bullet$}
\put (300,8){$\bullet$}
\put (325,33){$\bullet$}

\put (300,63){$0$}
\put (270,33){$1$}
\put (300,0){$2$}
\put (330,33){$3$}

\qbezier (277,35) (289.5,22.5) (302,10)
\qbezier (302,60) (302,35) (302,10)
\qbezier (302,60) (315.5,47.5) (327,35)

\end{picture}

\begin{center}
$
(X_{0,3},X_{1,2},X_{1,3})\hspace{7mm}(X_{2,3},X_{0,1},X_{0,2})\hspace{7mm}(X_{1,3},X_{2,3},X_{0,1})\hspace{7mm}(X_{0,2},X_{0,3},X_{1,2})$
\end{center}
\begin{center}
$
(X_{1,2},X_{0,3},X_{1,3})\hspace{7mm}(X_{0,1},X_{2,3},X_{0,2})\hspace{7mm}(X_{1,3},X_{0,1},X_{2,3})\hspace{7mm}(X_{0,2},X_{1,2},X_{0,3})$
\end{center}

We remark that the figures in the first column are same up to rotation.
Thus, we can see that there are 4 figures up to rotation.
In the next section, we consider the rotation of non-crossing spanning trees.

\end{ex}


\section{Action of the cyclic group}

Let $\Z / (n+1)\Z =\langle \s \rangle$ be the cyclic group of order $n+1$.
The group $\Z / (n+1)\Z$ acts on $C_{n+1}$ by

$$
\s (c(i,j))=\left\{
\begin{array}{ll}
c(i-1,j-1) & (i\not=0), \\
c(j-1,n) & (i=0).
\end{array}
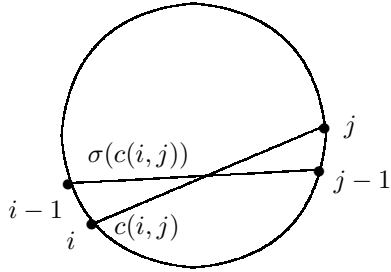
\right.
$$

\begin{center}
\begin{picture}(160,140)

\qbezier (20,65) (24,111) (70,115)
\qbezier (20,65) (24,19) (70,15)
\qbezier (70,115) (116,111) (120,65)
\qbezier (70,15) (116,19) (120,65)

\put (29,29){$\bullet$}
\put (117,65){$\bullet$}
\put (20,44){$\bullet$}
\put (115,49){$\bullet$}

\put (22,22){$i$}
\put (127,66){$j$}
\put (0,34){$i-1$}
\put (123,45){$j-1$}

\qbezier (33,32) (33,32) (118,68)
\qbezier (24,47) (24,47) (116,52)

\put (40,28){$c(i,j)$}
\put (30,55){$\s (c(i,j))$}

\end{picture}

\end{center}

The group $\Z / (n+1)\Z$ also acts on $\Gv$ via $\Phi$ and its action is given by

$$
\s (X_{i,j})=\left\{ \begin{array}{ll}
 X_{i-1,j-1} & (i\not=0)\\
 X_{j-1,n} & (i=0)
\end{array}
\right.
 \text{ for } X_{i,j} \in \Gv
$$

Thus, we see

$$
\s (X)=\left\{ \begin{array}{ll}
\tau X & {\rm if\ } X {\rm \ is\ not\ projective,}\\
\nu X & {\rm if\ } X {\rm \ is\ projective,}
\end{array}
\right.
$$
where $\nu$ is the Nakayama functor.

We set $\s (\e ) =(\s (E_1), \s (E_2), \cdots \s (E_n))$ for $\e=(E_1,E_2,\cdots ,E_n) \in \Exc$.
In general, we must check that $\s (\e )$ is a complete exceptional sequence by using the Auslander-Reiten duality.
But this is immediate from Theorem \ref{main} since $\Phi(\sigma(\epsilon))$ is clearly a non-crossing spanning tree.

\begin{ex}
If $n=4$, there are $11$ non-crossing spanning trees up to rotation.
For each figures, it has $5$ distinct rotations.
Thus, there are $55$ non-crossing spanning trees and $125$ exceptional sequences.


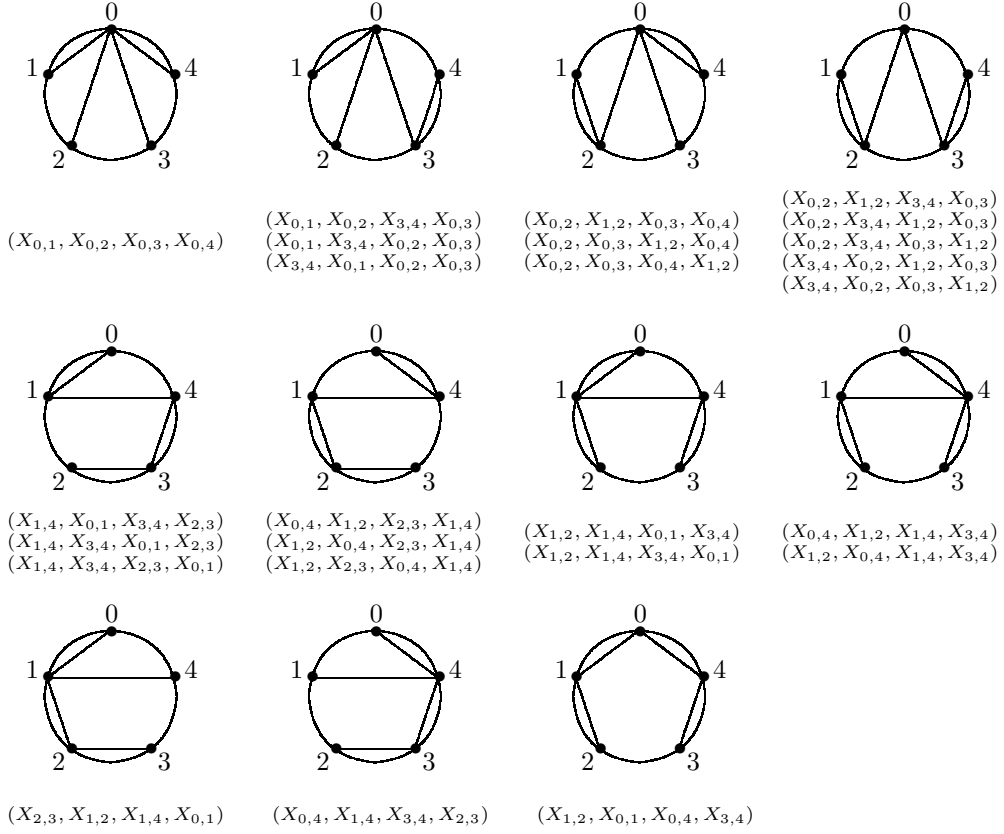
\begin{picture}(360,70)

\qbezier (7,35) (9,58) (32,60)
\qbezier (7,35) (9,12) (32,10)
\qbezier (32,60) (55,58) (57,35)
\qbezier (32,10) (55,12) (57,35)

\put (30,57){$\bullet$}
\put (6,40){$\bullet$}
\put (54,40){$\bullet$}
\put (15,13){$\bullet$}
\put (45,13){$\bullet$}

\put (30,63){$0$}
\put (0,42){$1$}
\put (10,7){$2$}
\put (50,7){$3$}
\put (60,42){$4$}

\qbezier (8,42) (20,51) (32,60)
\qbezier (17,15) (24.5,37.5) (32,60)
\qbezier (47,15) (39.5,37.5) (32,60)
\qbezier (56,42) (44,51) (32,60)

%

\qbezier (107,35) (109,58) (132,60)
\qbezier (107,35) (109,12) (132,10)
\qbezier (132,60) (155,58) (157,35)
\qbezier (132,10) (155,12) (157,35)

\put (130,57){$\bullet$}
\put (106,40){$\bullet$}
\put (154,40){$\bullet$}
\put (115,13){$\bullet$}
\put (145,13){$\bullet$}

\put (130,63){$0$}
\put (100,42){$1$}
\put (110,7){$2$}
\put (150,7){$3$}
\put (160,42){$4$}

\qbezier (108,42) (120,51) (132,60)
\qbezier (117,15) (124.5,37.5) (132,60)
\qbezier (147,15) (139.5,37.5) (132,60)
\qbezier (147,15) (152.5,28.5) (156,42)

%

\qbezier (207,35) (209,58) (232,60)
\qbezier (207,35) (209,12) (232,10)
\qbezier (232,60) (255,58) (257,35)
\qbezier (232,10) (255,12) (257,35)

\put (230,57){$\bullet$}
\put (206,40){$\bullet$}
\put (254,40){$\bullet$}
\put (215,13){$\bullet$}
\put (245,13){$\bullet$}

\put (230,63){$0$}
\put (200,42){$1$}
\put (210,7){$2$}
\put (250,7){$3$}
\put (260,42){$4$}

\qbezier (208,42) (212.5,28.5) (217,15)
\qbezier (217,15) (224.5,37.5) (232,60)
\qbezier (247,15) (239.5,37.5) (232,60)
\qbezier (256,42) (244,51) (232,60)

%

\qbezier (307,35) (309,58) (332,60)
\qbezier (307,35) (309,12) (332,10)
\qbezier (332,60) (355,58) (357,35)
\qbezier (332,10) (355,12) (357,35)

\put (330,57){$\bullet$}
\put (306,40){$\bullet$}
\put (354,40){$\bullet$}
\put (315,13){$\bullet$}
\put (345,13){$\bullet$}

\put (330,63){$0$}
\put (300,42){$1$}
\put (310,7){$2$}
\put (350,7){$3$}
\put (360,42){$4$}

\qbezier (308,42) (312.5,28.5) (317,15)
\qbezier (317,15) (324.5,37.5) (332,60)
\qbezier (347,15) (339.5,37.5) (332,60)
\qbezier (356,42) (351.5,28.5) (347,15)

\end{picture}

\noindent
{\scriptsize 
$
\begin{array}{l}
(X_{0,1},X_{0,2},X_{0,3},X_{0,4})
\end{array}
$\ \ \ 
$
\begin{array}{l}
(X_{0,1},X_{0,2},X_{3,4},X_{0,3})\\
(X_{0,1},X_{3,4},X_{0,2},X_{0,3})\\
(X_{3,4},X_{0,1},X_{0,2},X_{0,3})\\
\end{array}
$\ \ \ 
$
\begin{array}{l}
(X_{0,2},X_{1,2},X_{0,3},X_{0,4})\\
(X_{0,2},X_{0,3},X_{1,2},X_{0,4})\\
(X_{0,2},X_{0,3},X_{0,4},X_{1,2})\\
\end{array}
$\ \ \ 
$
\begin{array}{l}
(X_{0,2},X_{1,2},X_{3,4},X_{0,3})\\
(X_{0,2},X_{3,4},X_{1,2},X_{0,3})\\
(X_{0,2},X_{3,4},X_{0,3},X_{1,2})\\
(X_{3,4},X_{0,2},X_{1,2},X_{0,3})\\
(X_{3,4},X_{0,2},X_{0,3},X_{1,2})\\
\end{array}
$}

\vspace{10pt}


\begin{picture}(360,70)

\qbezier (7,35) (9,58) (32,60)
\qbezier (7,35) (9,12) (32,10)
\qbezier (32,60) (55,58) (57,35)
\qbezier (32,10) (55,12) (57,35)

\put (30,57){$\bullet$}
\put (6,40){$\bullet$}
\put (54,40){$\bullet$}
\put (15,13){$\bullet$}
\put (45,13){$\bullet$}

\put (30,63){$0$}
\put (0,42){$1$}
\put (10,7){$2$}
\put (50,7){$3$}
\put (60,42){$4$}

\qbezier (8,42) (20,51) (32,60)
\qbezier (8,42) (32,42) (56,42)
\qbezier (47,15) (32,15) (17,15)
\qbezier (56,42) (51.5,28.5) (47,15)


\qbezier (107,35) (109,58) (132,60)
\qbezier (107,35) (109,12) (132,10)
\qbezier (132,60) (155,58) (157,35)
\qbezier (132,10) (155,12) (157,35)

\put (130,57){$\bullet$}
\put (106,40){$\bullet$}
\put (154,40){$\bullet$}
\put (115,13){$\bullet$}
\put (145,13){$\bullet$}

\put (130,63){$0$}
\put (100,42){$1$}
\put (110,7){$2$}
\put (150,7){$3$}
\put (160,42){$4$}

\qbezier (156,42) (144,51) (132,60)
\qbezier (108,42) (132,42) (156,42)
\qbezier (108,42) (112.5,28.5) (117,15)
\qbezier (117,15) (132,15) (147,15)


\qbezier (207,35) (209,58) (232,60)
\qbezier (207,35) (209,12) (232,10)
\qbezier (232,60) (255,58) (257,35)
\qbezier (232,10) (255,12) (257,35)

\put (230,57){$\bullet$}
\put (206,40){$\bullet$}
\put (254,40){$\bullet$}
\put (215,13){$\bullet$}
\put (245,13){$\bullet$}

\put (230,63){$0$}
\put (200,42){$1$}
\put (210,7){$2$}
\put (250,7){$3$}
\put (260,42){$4$}

\qbezier (208,42) (220,51) (232,60)
\qbezier (208,42) (232,42) (256,42)
\qbezier (208,42) (212.5,28.5) (217,15)
\qbezier (256,42) (251.5,28.5) (247,15)


\qbezier (307,35) (309,58) (332,60)
\qbezier (307,35) (309,12) (332,10)
\qbezier (332,60) (355,58) (357,35)
\qbezier (332,10) (355,12) (357,35)

\put (330,57){$\bullet$}
\put (306,40){$\bullet$}
\put (354,40){$\bullet$}
\put (315,13){$\bullet$}
\put (345,13){$\bullet$}

\put (330,63){$0$}
\put (300,42){$1$}
\put (310,7){$2$}
\put (350,7){$3$}
\put (360,42){$4$}

\qbezier (356,42) (344,51) (332,60)
\qbezier (308,42) (332,42) (356,42)
\qbezier (308,42) (312.5,28.5) (317,15)
\qbezier (356,42) (351.5,28.5) (347,15)

\end{picture}

\noindent
{\scriptsize $
\begin{array}{l}
(X_{1,4},X_{0,1},X_{3,4},X_{2,3})\\
(X_{1,4},X_{3,4},X_{0,1},X_{2,3})\\
(X_{1,4},X_{3,4},X_{2,3},X_{0,1})\\
\end{array}
$\ \ \ 
$
\begin{array}{l}
(X_{0,4},X_{1,2},X_{2,3},X_{1,4})\\
(X_{1,2},X_{0,4},X_{2,3},X_{1,4})\\
(X_{1,2},X_{2,3},X_{0,4},X_{1,4})\\
\end{array}
$\ \ \ 
$
\begin{array}{l}
(X_{1,2},X_{1,4},X_{0,1},X_{3,4})\\
(X_{1,2},X_{1,4},X_{3,4},X_{0,1})\\
\end{array}
$\ \ \ 
$
\begin{array}{l}
(X_{0,4},X_{1,2},X_{1,4},X_{3,4})\\
(X_{1,2},X_{0,4},X_{1,4},X_{3,4})\\
\end{array}
$}

\vspace{10pt}


\begin{picture}(360,70)

\qbezier (7,35) (9,58) (32,60)
\qbezier (7,35) (9,12) (32,10)
\qbezier (32,60) (55,58) (57,35)
\qbezier (32,10) (55,12) (57,35)

\put (30,57){$\bullet$}
\put (6,40){$\bullet$}
\put (54,40){$\bullet$}
\put (15,13){$\bullet$}
\put (45,13){$\bullet$}

\put (30,63){$0$}
\put (0,42){$1$}
\put (10,7){$2$}
\put (50,7){$3$}
\put (60,42){$4$}

\qbezier (8,42) (20,51) (32,60)
\qbezier (8,42) (32,42) (56,42)
\qbezier (8,42) (12.5,28.5) (17,15)
\qbezier (17,15) (32,15) (47,15)


\qbezier (107,35) (109,58) (132,60)
\qbezier (107,35) (109,12) (132,10)
\qbezier (132,60) (155,58) (157,35)
\qbezier (132,10) (155,12) (157,35)

\put (130,57){$\bullet$}
\put (106,40){$\bullet$}
\put (154,40){$\bullet$}
\put (115,13){$\bullet$}
\put (145,13){$\bullet$}

\put (130,63){$0$}
\put (100,42){$1$}
\put (110,7){$2$}
\put (150,7){$3$}
\put (160,42){$4$}

\qbezier (156,42) (144,51) (132,60)
\qbezier (108,42) (132,42) (156,42)
\qbezier (156,42) (151.5,28.5) (147,15)
\qbezier (117,15) (132,15) (147,15)


\qbezier (207,35) (209,58) (232,60)
\qbezier (207,35) (209,12) (232,10)
\qbezier (232,60) (255,58) (257,35)
\qbezier (232,10) (255,12) (257,35)

\put (230,57){$\bullet$}
\put (206,40){$\bullet$}
\put (254,40){$\bullet$}
\put (215,13){$\bullet$}
\put (245,13){$\bullet$}

\put (230,63){$0$}
\put (200,42){$1$}
\put (210,7){$2$}
\put (250,7){$3$}
\put (260,42){$4$}

\qbezier (208,42) (212.5,28.5) (217,15)
\qbezier (208,42) (220,51) (232,60)
\qbezier (247,15) (251.5,28.5) (256,42)
\qbezier (256,42) (244,51) (232,60)

%
%
%
%

\end{picture}

\noindent
{\scriptsize 
$
\begin{array}{l}
(X_{2,3},X_{1,2},X_{1,4},X_{0,1})
\end{array}
$\ \ \ 
$
\begin{array}{l}
(X_{0,4},X_{1,4},X_{3,4},X_{2,3})
\end{array}
$\ \ \ 
$
\begin{array}{l}
(X_{1,2},X_{0,1},X_{0,4},X_{3,4})
\end{array}
$}

\end{ex}


\section{Mutation of exceptional sequences}

Let $(E,F)$ be an exceptional pair (in $\D$).
The {\it left mutation} $\l _EF$ of $F$ by $E$ is given by a distinguished triangle
$$ \bigoplus _{l \in \Z} \Hom (\S ^lE,F)\otimes _k \S ^lE \overset{can}{\to} F \to \l _EF,$$
where {\it can} denote the canonical morphism.
Dually, the {\it right mutation} $\r _FE$ of $E$ by $F$ is given by a distinguished triangle
$$ \r _FE \to E \overset{can}{\to} \bigoplus _{l \in \Z} \dual \Hom (E,\S ^lF)\otimes _k \S ^lF,$$
where $\dual$ is the $k$-dual functor $\Hom_k(-,k)$.

We denote by $B_n = \langle \s_1, \s_2, \cdots ,\s_{n-1}\ |\ \s_i\s_j=\s_j\s_i\ (j-i \geq 2),\ \s_i\s_{i+1}\s_i=\s_{i+1}\s_i\s_{i+1}\ (i=1,2,\cdots n-2) \rangle$ the braid group on $n$-strings and by $\Z^n=\bigoplus _{i=1}^n{\bf e}_i\Z$ a free abelian group of rank $n$.
It is known that $G:=B_n\ltimes \Z^n$ acts on $\Exc$ by
$$
\begin{array}{l}
\s_i (\e )= (E_1,E_2,\cdots ,E_{i-1},\l _{E_i}E_{i+1}, E_i,E_{i+2}, \cdots E_n),\\
{\bf e}_i (\e )= (E_1,E_2,\cdots ,E_{i-1},\S E_i,E_{i+1},\cdots E_n).
\end{array}
$$

Crawley-Boevey \cite{C} showed that this action is transitive, namely, for any $\e ,\e ' \in \Exc$, there exists $\rho \in G$ such that $\e = \rho \e '$.

In this section, we consider a mutation of chords which corresponds to exceptional pairs.

We extend $\Phi : \Gv =\{ \text{indecomposable $\R$-modules}\} \to C_{n+1}$ to $\Phi : \{ \text{indecomposable objects in }\D \} \to C_{n+1}$ by $\Phi (\S X)=\Phi (X)$.
By definition we have $\l _{\S ^iE}\S ^jF \cong \S ^j\l _EF$ and $\r _{\S ^jF}\S ^iE \cong \S ^i\r _FE$.
Thus we obtain $\Phi (\l _{\S ^iE}\S ^jF) = \Phi (\l _EF)$ and $\Phi (\r _{\S ^jF}\S ^iE) = \Phi (\r _FE)$.

The following lemma is straightforward.

\begin{lem}\label{5.1}
The following conditions are equivalent for an exceptional pair $(E,F)$:
\begin{enumerate}[1.]
\item $\l _EF \cong F$;
\item $\r _FE \cong E$;
\item $(F,E)$ is an exceptional pair;
\item $\bigoplus _{l \in \Z} \Hom (\S ^lE,F)=0$;
\item $\Phi (E)$ does not meet $\Phi (F)$.
\end{enumerate}
\end{lem}

From now on, we assume that $(E,F)$ is an exceptional pair but $(F,E)$ is not an exceptional pair.
Note that $\Phi (E)$ meets $\Phi (F)$ at some point $i$ by Lemma \ref{3.2} (3).

\begin{lem}\label{5.2}
Let $(E,F)$ be an exceptional pair with $\Phi (E)=c(i,j)$ and $\Phi (F)=c(i,l)$.
Then $\Phi (\l _EF) = \Phi (\r _FE) =c(j,l)$.
\end{lem}

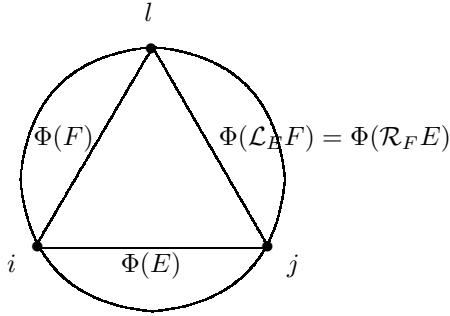
\begin{center}
\begin{picture}(140,140)

\qbezier (10,65) (14,111) (60,115)
\qbezier (10,65) (14,19) (60,15)
\qbezier (60,115) (106,111) (110,65)
\qbezier (60,15) (106,19) (110,65)

\put (57,112){$\bullet$}
\put (14,37){$\bullet$}
\put (101,37){$\bullet$}

\put (57,125){$l$}
\put (5,30){$i$}
\put (111,30){$j$}

\qbezier (16,40) (38,77.5) (60,115)
\qbezier (16,39) (60,39) (104,39)
\qbezier (60,115) (82,77) (104,39)

\put (48,30){$\Phi (E)$}
\put (15,78){$\Phi (F)$}
\put (85,78){$\Phi (\l _EF)=\Phi (\r _FE)$}

\end{picture}
\end{center}

\begin{pf}
There exist indecomposable $\R$-modules $X$ and $Y$, and integers $a$ and $b$ such that $E \cong \S ^aX$ and $F \cong \S ^bY$.
Since $\Hom (\S ^lE,F) \cong \Hom (\S ^l\S ^aX,\S ^bY) \cong \Hom (\S ^{l+a-b}X,Y)$ and 
$$
\Hom (\S ^{l+a-b}X,Y) \cong \left\{
\begin{array}{ll}
\Hom (X,Y) & (l+a-b=0), \\
\Ext (X,Y) & (l+a-b=-1), \\
0 & \text{otherwise},
\end{array}
\right.
$$
we can check that $\sum _{l \in \Z} \dim _k \Hom (\S ^lE,F) \leq 1$ by \cite[Lemma 8.1]{Hi} and by the Auslander-Reiten duality.
The assumption implies $\sum _{l \in \Z} \dim _k \Hom (\S ^lE,F) = 1$.
We may assume that $\Hom (E,F) \cong k$ and that $f$ is a generator of $\Hom (E,F)$.
Then, $\l _EF$ and $\r _FE$ are given by the following triangles:

$$
\begin{array}{ll}
E \overset{f}{\to} F \to \l _EF, \\
\r _FE \to E \overset{f}{\to} F.
\end{array}
$$
Therefore we get $\r _FE \cong \S ^{-1} \l _EF$ and $\Phi (\r _FE)=\Phi (\l _EF)$.

To see $\Phi (\l _EF)(=\Phi (\r _FE))=c(j,l)$, we assume $E$ and $F$ are $\R$-modules.

Note that if $i<l<j$, $l<j<i$ or $j<i<l$, then $(E,F)$ is not an exceptional pair by Lemma \ref{3.2}.

\begin{c1}
We assume $i<j<l$.
Then we have $E=X_{i,j}$ and $F=X_{i,l}$.
In this case, the generator of $\Hom (E,F)$ is a monomorphism and $\l _EF$ is obtained by a short exact sequence $0 \to E \to F \to \l _EF \to 0$.
Therefore we have $\l _EF =X_{j,l}$.
\end{c1}

\begin{c2}
We assume $j<l<i$.
Then we have $E=X_{j,i}$ and $F=X_{l,i}$.
In this case, the generator of $\Hom (E,F)$ is an epimorphism and $\r _FE$ is obtained by a short exact sequence $0 \to \r _FE \to E \to F \to 0$.
Therefore we have $\r _FE =X_{j,l}$.
\end{c2}

\begin{c3}
We assume $l<i<j$.
Then we have $E=X_{i,j}$ and $F=X_{l,i}$.
In this case, $\Ext (E,F) \cong k$ and $\l _FE$ is obtained by a short exact sequence $0 \to F \to \l _FE \to E \to 0$.
Therefore we have $\l _EF =X_{l,j}$.
\end{c3}

Thus we have $\Phi (\l _EF) = \Phi (\r _FE) = c(j,l)$.

\qed
\end{pf}

Now we give a definition of mutation of chords which corresponds to exceptional pairs.

\begin{dfn}For $c,c' \in C_{n+1}$, we define the {\it left mutation} $\l _cc'$ of $c'$ by $c$ and the {\it right mutation} $\r _{c'}c$ of $c$ by $c'$ as follows:
\begin{enumerate}[1.]
\item If $c$ does not meet $c'$, then we define $\l _cc'=c'$ and $\r_{c'}c=c.$
\item If $c=c(i,j)$ and $c'=c(i,l)$, then we define $\l _cc' = \r_{c'}c= c(j,l).$
\end{enumerate}
\end{dfn}

We can check the following Theorem by using Lemma \ref{5.1} and Lemma \ref{5.2}.

\begin{thm}
For an exceptional pair $(E,F)$, we have $\Phi (\l _EF) = \l _{\Phi (E)}\Phi (F)$ and $\Phi (\r _FE) = \r _{\Phi (F)}\Phi (E)$.
\end{thm}


\begin{ac}

The author is grateful Osamu Iyama and Frederic Chapoton for their comments and suggestions.
He thanks Hugh Thomas for pointing out that this paper is closely related to \cite{IT} and \cite{GY}.
He also thanks Ryo Takahashi and the referee for their careful reading.
\end{ac}


 
\ifx\undefined\bysame 
\newcommand{\bysame}{\leavevmode\hbox to3em{\hrulefill}\,} 
\fi


\begin{thebibliography}{1} 


\bibitem{A} T. Araya,
{\em Exceptional sequences over graded Cohen-Macaulay rings},
Math. J. Okayama Univ. {\bf 41} (1999), pp 81-102


\bibitem{ASS} I. Assem, D. Simson, and A. Skowronski, 
{\em Elements of the representation theory of associative algebras Vol. 1.
Techniques of representation theory},
 London Mathematical Society Student Texts, 65. Cambridge University Press, Cambridge, 2006.


\bibitem{ARS} M. Auslander, I. Reiten and S. O. Smalo,
{\em Representation theory of Artin algebras},
 Cambridge Studies in Advanced Mathematics, 36. Cambridge University Press, Cambridge, 1995.


\bibitem{BK} A. I. Bondal and  M. M. Kapranov,
{\em Representable functors, Serre functors, and reconstructions},
Izv. Akad. Nauk SSSR Ser. Mat. {\bf 53} (1989), no. 6, 1183--1205, 1337;
translation in Math. USSR-Izv. {\bf 35} (1990), no. 3, 519--541 


\bibitem{CCS} P. Caldero, F. Chapoton and  R. Schiffler,
{\em Quivers with relations arising from clusters ($A\sb n$ case)},
Trans. Amer. Math. Soc.  {\bf 358} (2006),  no. 3, 1347--1364.


\bibitem{C} W. Crawley-Boevey,
{\em Exceptional sequences of representations of quivers},
Proceedings of the Sixth International Conference on Representations of Algebras (Ottawa, ON, 1992),  7, 117--124, Carleton-Ottawa Math. Lecture Note Ser., 14, Carleton Univ., Ottawa, ON, 1992.


\bibitem{DP} S. Dulucq and J.-G. Penaud,
{\em Cordes, arbres et permutations},
Discrete Math. {\bf 117} (1993), 89--105.


\bibitem{FN} P. Flajolet and M. Noy,
{\em Analytic combinatorics of non-crossing configurations},
Discrete Mathematics {\bf 204} (1999), 203--229

\bibitem{G} P. Gabriel,
{\em Auslander-Reiten sequences and representation-finite algebras},
Representation theory, I (Proc. Workshop, Carleton Univ., Ottawa, Ont., 1979),
pp. 1--71, Lecture Notes in Math., 831, Springer, Berlin, 1980.


\bibitem{GR} A. L. Gorodentsev and A. N. Rudakov,
{\em Exceptional vector bundles on projective spaces},
Duke. math. J. {\bf 54} (1987), 115--130.


\bibitem{GY} I. Goulden and A. Yong,
{\em Tree-like properties of cycle factorizations},
 J. Combinatorial Theory (A) {\bf 98}, 2002, 106--117. 


\bibitem{H} D. Happel,
{\em Triangulated categories in the representation theory of finite-dimensional algebras},
 London Mathematical Society Lecture Note Series, 119. Cambridge University Press, Cambridge, 1988


\bibitem{Hi} L. Hille,
{\em On the volume of a tilting module},
 Abh. Math. Sem. Univ. Hamburg {\bf 76} (2006), 261--277.


\bibitem{IT} C. Ingalls and H. Thomas,
{\em Noncrossing partitions and representations of quivers},
Preprint (2008), \texttt{http://arxiv.org/abs/math/0612219v4}.


\bibitem{I} O. Iyama,
{\em Higher-dimensional Auslander-Reiten theory on maximal orthogonal subcategories},
Adv. Math.  {\bf 210}  (2007),  no. 1, 22--50.


\bibitem{M} H. Meltzer,
{\em Exceptional vector bundles, tilting sheaves and tilting complexes for weighted projective lines},
Mem. Amer. Math. Soc. {\bf 171} (2004), no. 808


\bibitem{Ri} C. Riedtmann,
{\em Representation-finite self-injective algebras of class $A\sb{n}$},
Representation theory, II (Proc. Second Internat. Conf., Carleton Univ., Ottawa, Ont., 1979),  pp. 449--520, Lecture Notes in Math., 832, Springer, Berlin, 1980.

\bibitem{Rin} C. M. Ringel,
{\em The braid group action on the set of exceptional sequences of a hereditary Artin algebra},
 Abelian group theory and related topics (Oberwolfach, 1993), 339--352, Contemp. Math., {\bf 171}, Amer. Math. Soc., Providence, RI, 1994.

\bibitem{R}  A. N. Rudakov,
{\em Exceptional collections, mutations and helices},
Helices and vector bundles,  1--6, London Math. Soc. Lecture Note
Ser., 148, Cambridge Univ. Press, Cambridge, 1990.


\bibitem{S} U. Seidel,
{\em Exceptional sequences for quivers of Dynkin type},
Comm. Algebra {\bf 29} (3) (2001), 1373--1386.


\bibitem{W} A. Wiedemann,
{\em Die Auslander-Reiten K\"oher der gitterendlichen Gorensteinordnungen},
Bayreuth. Math. Schr.  No. {\bf 23} (1987), 1--134.

\end{thebibliography}
\end{document}